\documentclass [A4paper]{article}
\usepackage{amstext,amsthm,latexsym,amsmath, amsfonts,amssymb, amscd, bezier}
\usepackage[latin1]{inputenc}
\usepackage[mathscr]{eucal}
\usepackage{graphicx}
\usepackage{epsfig}
\usepackage{color}
\usepackage[normalem]{ulem}
\usepackage[all,cmtip]{xy}

\newcommand{\dis}{\displaystyle}
\newcommand{\C}{\mathbb{C}}

\newcommand{\Coker}{\mathop{\mathrm{Coker}}}
\def\cO{\mathcal O}

\newtheorem{proposition}{Proposition}[section]
  \newtheorem{theorem}[proposition]{Theorem}
  
  \newtheorem{example}[proposition]{Example}
  \newtheorem{defi}[proposition]{Definition}
  \newtheorem{remark}[proposition]{Remark}

\begin{document}

\title{An algorithm to compute a presentation of pushforward modules}
\author{M. E. Hernandes, A. J. Miranda, G. Peñafort-Sanchis}

\maketitle

\begin{abstract}
We describe an algorithm to compute a presentation of the
pushforward module $f_*\cO_{\mathcal{X}}$ for a finite map germ
$f\colon \mathcal{X}\to (\C^{n+1},0)$, where $\mathcal{X}$ is
Cohen-Macaulay of dimension $n$. The algorithm is an improvement of
a method by Mond and Pellikaan. We give applications to problems in
singularity theory, computed by means of an implementation in the software {\sc Singular}.
\end{abstract}

\begin{center} Mathematics Subject Classification. Primary 58K05; Secondary 32S10, 14Q99, 13C10 .\end{center}

\begin{center} Keywords: presentation matrices, Fitting ideals, multiple-point schemes. \end{center}

\section{Introduction}

Let $M$ be a module over a commutative unitary ring  $R$. A
\emph{presentation} of $M$  is an exact sequence
$$R^p \stackrel{\lambda}{\longrightarrow} R^q \stackrel{\psi}\longrightarrow M \longrightarrow 0.$$
If $M$ admits a presentation, then we say that $M$ is a \emph{finitely presented} module, and any matrix $\Lambda$ associated to $\lambda$ is called a
\emph{presentation matrix} of $M$.  It is well known that any finitely generated module over a Noetherian ring is finitely presented (see \cite{GP}).

 In this work, we introduce an algorithm to compute presentation matrices in the particular case where $M$ is the pushforward of the ring of holomorphic functions in the source of a finite map $f\colon \mathcal X\to (\C^{n+1},0)$, and $\mathcal X$ is a Cohen-Macaulay space of dimension $n$. The algorithm is based on a method by Mond and Pellikaan \cite{MP}, but introduces an improvement which allows to circumvent certain problems, concerning the limitation to polynomial inputs and outputs of commutative algebra
systems, such as {\sc Singular} \cite{singular}. As we will see, this improvement also makes the algorithm more efficient from a computational point of view. The reader can find in \cite{algoritmo} a {\sc Singular} library containing an
implementation of the algorithm.

 In the last section, we illustrate some applications of presentation matrices that can be computed by means of our algorithm. We show the computation of target and source multiple-point schemes for map germs $f\colon
(\C^n,0)\to(\C^{n+1},0)$, discriminants and certain topological invariants of maps ---leading, for example, to the answer of a question, due to Gaffney and Mond, about the
topological classification for corank $2$ map germs from
$\mathbb{C}^2$ to $\mathbb{C}^2$. All these applications are based on \emph{Fitting ideals}, which play a crucial role in the theory of singularities of map germs and  in
enumerative geometry (see for instance \cite{KLU}).

 If $M$ is a finitely presented module as above, then its $k$-th Fitting ideal is given by
\[{  F}_k(M)=\left \{
\begin{array}{ccl}
0 & \mbox{if} & k<0; \\ \langle \mbox{minors
 of order}\ q-k\ \mbox{of}\
\Lambda \rangle & \mbox{if} & 0\leq k<\min(q,p); \\
R & \mbox{if} & \min(q,p)\leq k.
\end{array}
\right .\]
Fitting ideals are invariant under module isomorphisms and they do
not depend on the chosen presentation of $M$ (see \cite{Lo}).

Let $f\colon \mathcal X\to \mathcal Y$ be a finite map germ and $\cO_{\mathcal X}, \cO_{\mathcal Y}$ the rings of holomorphic functions of $\mathcal X$ and $\mathcal Y$. The \emph{pushforward module}
$f_*\cO_{\mathcal X}$ is just $\cO_{\mathcal X}$, regarded as an $\cO_{\mathcal Y}$-module via $f$. Finiteness of $f$ implies that $f_*\cO_{\mathcal
X_1}$ is finite, and hence finitely presented. For simplicity, we write the corresponding Fitting ideals as $F_{k}(f)=
F_{k}(f_*\cO_{\mathcal X})$.  As shown by
Mond and Pellikaan \cite{MP}, the $k$th Fitting ideal of
$f_*\cO_{\mathcal X}$ defines the $(k-1)$th multiple point space
of $f$ in $\mathcal Y$, which we write as $M_k(f)=V(
F_{k-1}(f))$.

It is worth saying that the {\sc Singular} implementation of the algorithm has already been used in the works of other authors: In \cite{BOT}, Nu\~no Ballesteros, Or\'efice Okamoto and Tomazella use it to compute discriminant curves of map germs from
a complete intersection surface to the plane. Oset Sinha, Ruas and Wik Atique have used the algorithm to compute the image of a stable map germ of corank 2 from $\C^8$ to $\C^9$, which plays an important role in their work on the extra-nice dimensions \cite{ORW}.  Recently, O. N. Silva has used the algorithm to compute source double points of certain maps (as in Section \ref{secSourceDoublePoints}) \cite{S}, obtaining the first known counter-example to a conjecture by M. A. Ruas, on the equivalence between Whitney equisingularity and Topological triviality. For convenience to the reader, we describe Mond-Pellikaan's original method to obtain presentation matrices:

\subsection{Mond-Pellikaan algorithm}
Let $\mathcal X$ be an $n$-dimensional germ of Cohen Macaulay space, and let \[f\colon \mathcal X\to (\C^{n+1},0)\] satisfy the following extra condition: If we let \[\tilde f\colon \mathcal X\to (\C^{n},0),\] be the germ obtained by composing $f$ with the projection $(\C^{n+1},0)\to (\C^n,0)$ which forgets the last coordinate, then $\tilde f$ is a finite map germ.  Since our interest is to obtain a computer implementation, we are going to switch from the holomorphic setting of Mond and Pellikaan to the rational setting software like {\sc Singular} can handle. This is mostly a matter of language, and the results of Mond and Pellikaan apply here exactly in the same way.

Let $A=\C[X,Y]_{\langle X,Y\rangle}$ be the localization at the maximal ideal at the origin of the ring
of polynomials in the $n+1$ variables $X=X_1,\dots,X_n$ and $Y$. We
denote $\tilde A=\C[X]_{\langle X\rangle}$ and $B=(\C[x]/I)_{\langle
x\rangle}$, with variables $x=x_1,\dots,x_\ell$, and assume that $B$
is a Cohen-Macaulay ring of dimension $n$. Let
\[\phi\colon A\to B\] be a morphism of local rings given by
$X_i\mapsto f_i, i=1,\dots,n$ and $Y\mapsto f_{n+1}$, for some
polynomials  $f_j\in \C[x]$. Write
 \[\tilde \phi\colon \tilde A\to B\]
for the restricted morphism, and assume that $B$  is minimally generated by $g_1,\dots,g_h$ as an $\tilde A$-module.
Since $B$ is generated
 by $g_1,\dots,g_h$, there exist $\alpha_{ij}\in\tilde A$, $1\leq i,j\leq
h$, satisfying the equations
\begin{align}\label{equation}
 Y g_i=\sum_{j=1}^h \alpha_{ij}
g_j,\text{ for every }1\leq i\leq h.
\end{align}
Let $\lambda:A^h\to A^h
$ be given by multiplication by the matrix $\Lambda$ whose entries are
\[\Lambda_{ij}=\alpha_{ij}-\delta_{ij}Y,\]
where $\delta_{ij}$ stands for the Kronecker delta function.

If $\psi\colon A^h\to B$ is the epimorphism given by $e_i\mapsto
g_i$, where $e_i\in A^h$ is the element whose only non-zero entry is
$1$ in the $i$th position, then the inclusion $\mbox{Im } \lambda\subseteq \mbox{Ker } \psi$ follows from  (\ref{equation}). Mond and
Pellikaan show  that indeed the sequence
\begin{equation}\label{exactsequence} A^h \stackrel{\lambda}{\longrightarrow}
 A^h \stackrel{\psi}\longrightarrow B \longrightarrow 0\end{equation} is exact \cite{MP}. Therefore, the matrix
$$\Lambda=\left(\begin{array}{cccc}\alpha_{11}-Y & \alpha_{12} & \cdots & \alpha_{1h} \\ \alpha_{21} & \alpha_{22}-Y & \cdots & \alpha_{2h} \\ \vdots & \vdots & \ddots &  \vdots \\ \alpha_{h1} & \alpha_{h2} &\cdots  & \alpha_{hh}-Y\end{array}\right)$$
is a presentation matrix for $B$.
\begin{example}{\rm Let $f:(\mathbb C^2,0)\to(\mathbb C^3,0)$ be the cross-cap map (Figure 1), given by $$(x,y)\mapsto(x,y^2,xy).$$
Let $A=\C[X_1,X_2,Y]_{\langle X_1,X_2,Y\rangle}$, $\tilde
A=\C[X_1,X_2]_{\langle X_1,X_2\rangle}$ and $B=\C[x,y]_{\langle
x,y\rangle}$ and consider the ring homomorphism $A\to B$ given by
$X_1\mapsto x,X_2\mapsto y^2$ and $Y\mapsto xy$.

Observe that $B$ is an $\tilde A$-module minimally
generated by $g_1=1$ and $g_2=y$, and thus we can apply
Mond-Pellikaan algorithm. From the equalities $$Y\cdot g_1=xy\cdot g_1=x\cdot g_2=X_1\cdot
g_2\ \ \ \mbox{and}\ \ \ Y\cdot g_2=xy\cdot g_2=xy^2\cdot
g_1=X_1X_2\cdot g_1$$ we obtain the presentation matrix
$$\Lambda=\left(\begin{array}{cc}-Y & X_1 \\X_1X_2 & -Y\end{array}\right).$$
The matrix $\Lambda$ can be used to compute the Fitting ideals of $f_*\cO_2$, which determined the multiple-point scheme in the target of $f$.  The image of $f$ and the
double point space are, respectively:
 \[M_1(f)=V(Y^2-X_1^2X_2)\ \ \ \mbox{and}\ \
\ M_2(f)=V(X_1, Y).\] }
\begin{figure}
\begin{center}
\includegraphics[scale=0.45]{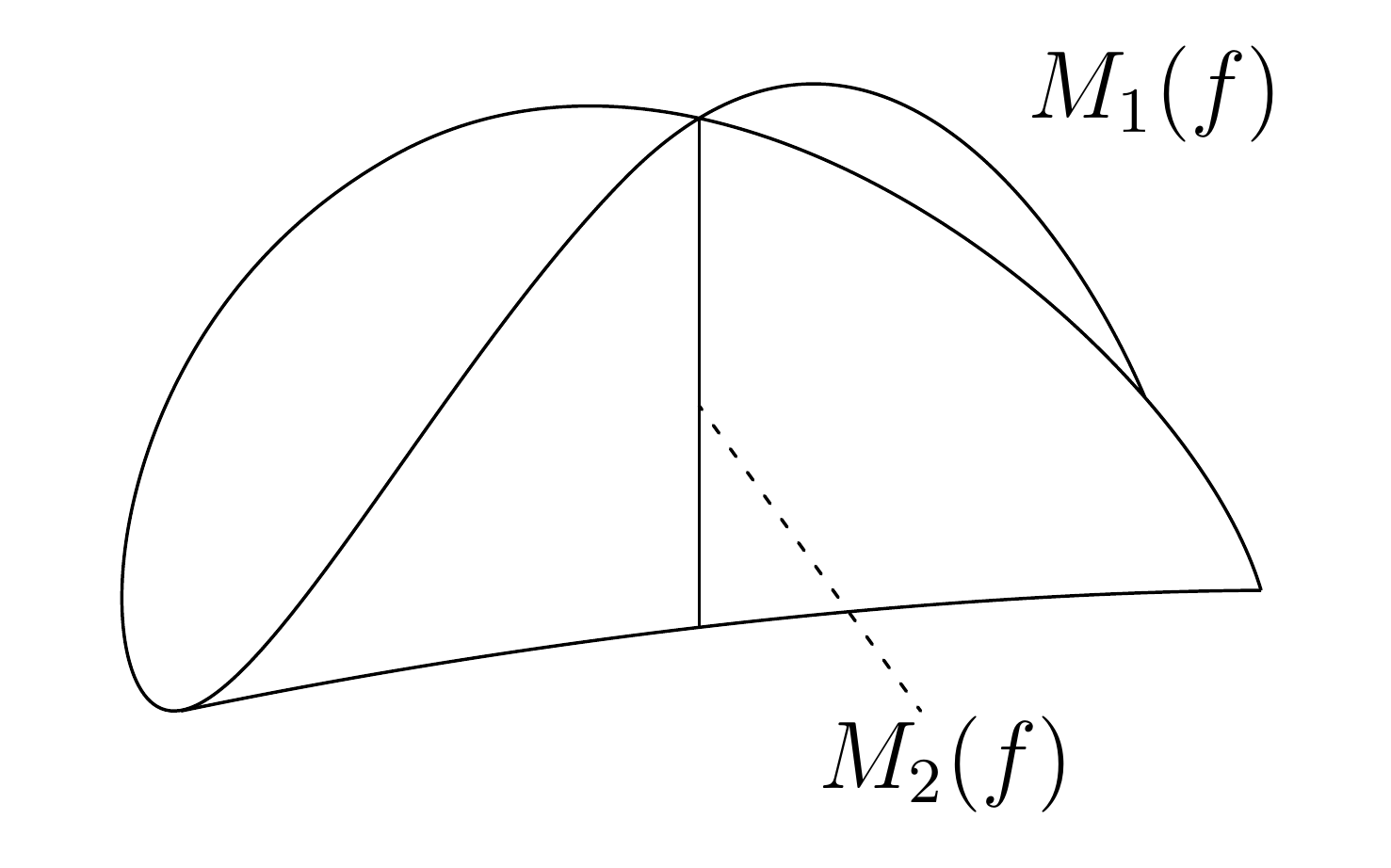}
\caption{ Image of a cross-cap.}
\end{center}
\end{figure}
\end{example}

\section{Polynomial presentation matrices}
Commutative algebra software such as {\sc Singular} only admit polynomial inputs and outputs. In this section we deal with the problem of how to find presentation matrices whose entries are polynomial. We start with the following trivial remark:
\begin{remark}\label{remPolyPresMat}
The elements $\alpha_{ij}\in \tilde A$ in (\ref{equation}) are fractions
$\alpha_{ij}=a_{ij}/b_{ij}$, for some polynomials $a_{ij},b_{ij}\in
\C[X]$   and $b_{ij}(0)\neq 0$. Multiplying the $i$th row by the
least common multiple of the elements $b_{ij},\ j=1,\dots, h$, we
obtain another presentation matrix, whose entries are polynomial.
\end{remark}

The previous remark guarantees, given a minimal collection of generators $g_i$, the
existence of a polynomial presentation matrix  of the following
form:

\begin{defi}
Given $g_1,\dots, g_h\in B$, an
\emph{$\bf{MP}$-matrix} (for $g_i$) is a matrix
$$\Lambda=\left(\begin{array}{cccc}\beta_{11}-u_{1}Y & \beta_{12} & \cdots & \beta_{1h} \\ \beta_{21} & \beta_{22}-u_{2}Y & \cdots & \beta_{2h} \\ \vdots & \vdots & \ddots &  \vdots \\ \beta_{h1} & \beta_{h2} &\cdots  & \beta_{hh}-u_{h}Y\end{array}\right),$$
with $\beta_{ij},u_{j}\in \C[X], u_{j}(0)=1$, such that the equations (\ref{equation}) hold.
\end{defi}

With the previous notations, let $g_1,\dots,g_h$ be a minimal system of generators of $B$ as an $\tilde A$-module. It follows from Remark \ref{remPolyPresMat} that $B$ admits an {\bf MP}-matrix $\Lambda$ for $g_1,\dots,g_h$ as a presentation matrix. The matrix is given by some polynomials $\beta_{ij},u_{j}\in \C[X], u_j(0)=1$, satisfying the conditions
\[\phi (u_jY)g_j\equiv\sum_{i=1}^h\phi(\beta_{ij})g_i \mod I,\text { for all }1\leq j\leq h.\]
To find the $j$th row of such $\Lambda$, one fixes polynomials up to
some degree $d$:$$\beta_{ij}=\sum_{\vert \alpha\vert\leq d}
a_{i,\alpha}X^\alpha \text{ for }i=1,\dots,h,$$
$$u_{j,\alpha}=1+\sum_{1\leq\vert \alpha\vert\leq d} b_\alpha
X^\alpha,$$ and tries to find $a_{i,\alpha},b_\alpha\in\mathbb{C}$,
such that the polynomial
\[P_d(x)=\phi
(u_jY)g_j-\sum_{i=1}^h\phi(\beta_{ij})g_i\]
 reduces to $0$ modulo
$I$. This is a linear system on $a_{i,\alpha},b_\alpha$, prescribed
by the vanishing of the coefficient of each $x^{\alpha}$ in the
reduction of $P_d$ modulo $I$. If that is not possible for the
degree $d$, then one increases $d$ and starts all over again.

Due to the reduction process, and to the clearing of denominators in Remark \ref{remPolyPresMat},
there is no obvious way to estimate the degrees of the entries in an
{\bf MP}-matrix $\Lambda$ in terms of the degrees of $f_i$ and the
generators of $I$. Unfortunately, the usage of reductions, and the
increasing number of parameters $a_{i,\alpha}, b_\alpha$ involved,
make the complexity of the procedure explained above grow very
rapidly as $d$ increases. In order to keep the degree $d$ as low as
possible, it seems a good idea to consider a class of matrices
bigger than the set of {\bf MP}-matrices. The following example
illustrates this situation.

\begin{example}\label{exHighLowDegree}{\rm
Let $\mathcal{X}=V(I)\subseteq \C^3$, with $I=\langle
z-x^ky\rangle$, and let $f\colon \mathcal{X}\to \C^3$ be given by
$$(x,y,z)\mapsto (x,y^2+xz,z).$$
In our usual setting $A=\mathbb{C}[X_1,X_2,Y]_{\langle
X_1,X_2,Y\rangle}, \widetilde{A}=\mathbb{C}[X_1,X_2]_{\langle
X_1,X_2\rangle}, B=\left (\frac{\mathbb{C}[x,y,z]}{\langle
y-x^ky\rangle}\right )_{\langle x,y,z\rangle}$, $\phi$ is given by
$X_1\mapsto x, X_2\mapsto y^2+xz$ and $Y\mapsto z$, and the
pushforward module $B$ is minimally generated by $g_1=1$ and $g_2=y$
as an $\tilde A$-module. It is easy to check that the matrix
$$\Lambda=\left(\begin{array}{cc}Y & -X_1^k \\-X_1^kX_2 & Y+X_1^{2k+1}\end{array}\right)$$
is an {\bf MP}-matrix for $g_1,g_2$. However, the matrix
$$\Lambda'=\left(\begin{array}{cc}Y & -X_1^k \\-X_1^k(X_1Y-X_2) & Y\end{array}\right)$$
is also a presentation matrix (by Theorem \ref{teorema} below). The
procedure explained before will have to consider polynomials up-to
degree $2k+1$ in order to find $\Lambda$, but a more flexible
version, allowing matrices such as $\Lambda'$, will stop at degree
$k+2$.}
\end{example}

\begin{defi}
 With the previous notations, given $g_1,\dots, g_h\in B$, an \emph{$\bf{HMP}$-matrix (for $g_i$)} is a matrix $\Lambda$ with polynomial entries in $\C[X,Y]$, satisfying the following conditions:

 \begin{itemize}
\item[{\bf C1.}]$\displaystyle\sum_{j=1}^{h} \phi(\Lambda_{ij})g_j\equiv 0 \mod I,$ for $ i=1,\dots,h$.
\item[{\bf C2.}]$\Lambda_{ij}(0,Y)-\Lambda_{ij}(0,0)=
\left \{ \begin{array}{l} Y\cdot u_i(Y); \mbox{for some polynomial}\
u(Y)\in\C[Y]\ \mbox{satisfying}\ u(0)\neq 0\ \mbox{if}\ i=j; \\
0\ \mbox{if}\ i\neq j.\end{array}\right .$
\end{itemize}
\end{defi}

\begin{remark}
Every {\bf MP}-matrix is an {\bf HMP}-matrix. On the
other hand, the {\bf HMP}-matrix $\Lambda'$ in Example \ref{exHighLowDegree}
is not an {\bf MP}-matrix.
\end{remark}

\begin{theorem}\label{teorema} If $\Lambda$ is an {\bf HMP}-matrix for a minimal set of generators $g_1,\dots,g_h$ of $B$ as an $\tilde A$-module, then $\Lambda$ is a presentation matrix for $B$ as an $A$-module.
\begin{proof}
The proof is similar to the one for Mond-Pellikaan's algorithm. Take the sequence of $A$-modules
\[A^h\stackrel{\Lambda}{\longrightarrow}A^h\stackrel{\psi}{\longrightarrow} B\longrightarrow 0,\]
where $\psi$ is determined by $e_i\mapsto g_i$, and $e_i$ the $i$th canonical vector in $A^h$. Condition {\bf C1} implies $\mbox{Im }
\Lambda\subseteq \mbox{Ker } \psi$, and $\psi$ is an epimorphism, so it suffices to show that $\Coker \Lambda= B$.

Let $B'=(\C[x,t]/I)_{\langle x,t\rangle}$ and let $\phi'\colon A\to B'$ be  given by
$X_i\mapsto f_i, i=1,\dots,n$ and $Y\mapsto f_{n+1}+t$. From the fact that $B$ is minimally generated
 by $g_1,\dots,g_h$ as an $\tilde A$-module, it follows that $B'$ is a free $A$-module minimally generated by $g_1,\dots,g_h$. Let
 $\eta \colon A^{h}\rightarrow
B'$ be the isomorphism given by
\[ (a_1,\ldots
,a_h)\mapsto\sum_{i=1}^{h}\phi'(a_i)g_i\]
and let $\varphi\colon B'\to B'$ be the morphism defined by extending  \[g_j\mapsto \sum_{i=1}^{h}\phi'(\Lambda_{ij})\cdot g_i,\text{ for }j=1,\ldots ,h,\] $A$-linearly, so that the diagram
\[\xymatrix{
A^{h}     \ar[d]^{\eta}                 \ar[r]^{\Lambda} &  A^{h} \ar[d]^{\eta} \\
            B' \ar[r]^{\varphi} &  B'  }\]
            commutes. We will show that $\Coker \varphi=B$.

The morphism $\varphi$ extends  $A$-linearly the assignations $g_j\mapsto \sum_{i=1}^{h}\phi'(\Lambda_{ij}) g_i$, for $j=1,\ldots ,h$. Consider the expansion
\[\begin{array}{lcl}
\dis\sum_{i=1}^{h}\phi'(\Lambda_{ij})\cdot g_i& = &\dis\sum_{i=1}^{h}(\Lambda_{ij}(f_1, \ldots, f_{n+1}+t))\cdot g_i\\
                                   & = &\dis\sum_{i=1}^{h}(\Lambda_{ij}(f_1, \ldots,  f_{n+1})) \cdot g_i+\dis\sum_{i=1}^{h}(\dis\sum_{k=1}^{\infty}\frac{1}{k!}
                                       \frac{\partial^k \Lambda_{ij}}{\partial Y^k}(f_1, \ldots,  f_{n+1})t^k)\cdot g_i\\
                                   & = &t  \dis\sum_{i=1}^{h}(\dis\sum_{k=1}^{\infty}\frac{1}{k!}
                                       \phi(\frac{\partial^k \Lambda_{ij}}{\partial Y^k})t^{k-1}) \cdot g_i.
\end{array}
\]
It follows that $\varphi$ splits as the composition $B'\stackrel{\cdot t}{\longrightarrow} B'\stackrel{\psi}{\longrightarrow}B'$, where the first morphism is multiplication by $t$ and $\psi$ is obtained by extending $g_j\mapsto g_j'=\sum_i R_{ij}g_i$, with
\[R_{ij}=\dis\sum_{k=1}^{\infty}\frac{1}{k!}\phi(\frac{\partial^k
\Lambda_{ij}}{\partial Y^k}) t^{k-1}.\]
It suffices to show that $\psi$ is an $A$-module isomorphism, that is, that $g_1',\dots,g_h'$ is a system of generators of $B'$. This is equivalent to show that the collection of the classes of $g_1',\dots,g_h'$ is a $\C$-basis $B'/\mathfrak m B'$, where $\mathfrak m$ is the maximal ideal in $A$. If $a\in A$ is divisible by some $X_i$, then $\phi(a)\in \mathfrak m B\subset\mathfrak m B'$. Therefore, condition {\bf C2} implies that the classes of the non-diagonal coefficients $R_{ij}, i\neq j$, are all zero, and the classes of $R_{ii}$ are all non-zero. This implies that the collection of classes of $g_1',\dots,g_h'$ is a basis, as desired.
\end{proof}
    \end{theorem}

\begin{remark}
It is clear that Theorem \ref{teorema} works for holomorphic maps as
well. With our usual assumptions, if a holomorphic map $f\colon
\mathcal X\to \C^{n+1}$ has polynomial coordinate functions  and
$g_1,\dots,g_h$ are polynomial minimal set of generators of $\tilde
f_*\cO_X$, then any {\bf HMP}-matrix for $g_1,\dots,g_h$ is a
presentation matrix for $f_*\cO_X$.
\end{remark}

\section{Algorithm and applications}

In this section we describe an algorithm to obtain a matrix
$\Lambda$ satisfying {\bf C1} and {\bf C2}, and we give some
applications. An implementation of this algorithm in {\sc Singular}
can be found in \cite{algoritmo}.

We recall previous assumptions and notations: We use variables $X=X_1,\dots,X_n$ and $Y$, and variables $x=x_1,\ldots ,x_{\ell}$. We write $A=\C[X,Y]$ and $\tilde A=\C[X]$. $I$ is an ideal in
$\mathbb C[x]_{\langle x\rangle}$, such that $B=\C[x]_{\langle x\rangle}/I$ is a Cohen Macaulay $n$-dimensional ring.
We have a morphism of rings $\phi\colon A\to B$, such that $B$ is finitely generated as $\tilde A$-module. It is well known that a minimal set of
generators $\{g_1,\dots, g_h\}$ for $B$ as $\tilde A$-module can be obtained
as representatives of the elements of a basis of the $\C$-vector
space $B/\mathfrak m B$, where $\mathfrak m$ is
the ideal maximal ideal $\langle X_1,\ldots ,X_n\rangle$ in $\tilde A$. We assume that such a basis can be
computed by an internal procedure of the software used for
implementation. We also assume the software to be able to perform
Groebner basis computations, in particular the reduction of an ideal
with respect to another one (see \cite{GP}). In {\sc Singular}, this
operations can be computed by using the instructions \texttt{kbase}
and \texttt{reduce}, respectively.

The outline of our algorithm is as follows:

\begin{center}
\begin{tabular}{|l|}
\hline
{\sc Inputs:}  $f$ and $I$. \\

Compute a $\mathbb{C}$-basis $\{g_1,g_2,\ldots ,g_h\}$ of $B/\mathfrak m B$. \\

{\sc For} $i=1,\ldots ,h$ do:\\

\hspace{0.5cm} {\sc Define} $w:=1$ and $k:=0$; \\

\hspace{0.5cm} {\sc While} $w\neq 0$ {\sc do}: \\

\hspace{1cm} $k:=k+1;$ \\

\hspace{1cm} Consider  $v_{i1},\ldots ,v_{ih}$, where\\

\hspace{1.7cm} $v_{ij}=\sum_{\vert \alpha \vert\leq k}a_{ij}^\alpha X^\alpha$ is a \emph{generic polynomial} satisfying {\bf (P)} (see Remark \ref{remark});\\

\hspace{1cm} Compute the reduction $P(a_{ij}^\alpha,x)$ of $\sum_{j=1}^{h}\phi(v_{ij})g_j$ modulo $I$;\\

\hspace{1cm} {\sc If} there exists $\tilde a_{ij}^\alpha\in \C$, such that $P(\tilde a_{ij}^\alpha,x)=0$ \\

\hspace{1.7cm} {\sc then} $\lambda_{ij}:=\sum_{\vert \alpha \vert\leq k}\tilde a_{ij}^\alpha X^\alpha$ and $w:=0$; \\

{\sc Output:} Matrix presentation $\Lambda=(\lambda_{ij})$. \\

\hline
\end{tabular}
\end{center}
\begin{center}{\sc Algorithm to compute an HMP-matrix.}\end{center}

\begin{remark}\label{remark}{\rm
By a generic polynomial we mean that the coefficients $a_{ij}^\alpha$ are parameters in the base ring. Property {\bf (P)} is as follows:
\begin{itemize}
\item $v_{ii}(0,\ldots ,0,Y)\equiv Y\mod\ (Y^2)$;
\item $v_{ij}(0,\ldots ,0,Y)\in\mathbb{C}$, for all $j\neq i$.
\end{itemize}
}
\end{remark}

By construction, condition {\bf (P)} implies that
$\Lambda$ satisfies condition {\bf C2}, and the fact that the
reduction $P(\tilde a_{ij}^\alpha,X)$ vanishes ensures that {\bf C1}
holds. Note that the degree of the generic polynomials $v_{ij}$
grows with the ``while'' loop, and the algorithm runs over all the
matrices considered in Remark  \ref{remPolyPresMat}. Since we assumed that $B$ is a finitely generated $\tilde A$-module, the algorithm terminates.

In the remaining subsections we show some applications of the
implementation \cite{algoritmo} of the above algorithm. All
computations were done using a computer with $2.8$ $Ghz$ Intel Core
$I7$ processor and $8$ $Gb$ of RAM memory. We use standard
singularity theory notation for which the reader can find the
details in the references.

\subsection{Topological invariants for  $f: (\C^2,0) \longrightarrow (\C^2,0)$}

Let $f: (\mathbb C^2,0) \longrightarrow (\mathbb C^2,0)$ be the map
germ given by $$(x,y)\longmapsto (xy,x^4+y^{37}+x^2y^{23}).$$ The
singular set is $\Sigma(f)=V(4x^4-21x^2y^{23}-37y^{37})$ and the
presentation matrix $\Lambda$ for $f_*\mathcal{O}_{\Sigma(f)}$ is a
$41 \times 41$ matrix, too big to be written here. Using our
implementation, the total time to obtain such matrix was about $95$
seconds.

An important topological invariant of $f$ is the Milnor number of
the discriminant $\mu(\Delta(f))$. Using {\sc Singular} software (in
the same computer), we can not obtain $\mu(\Delta(f))$ due to lack
of memory. But, in \cite{GM1} we find that
$\mu(\Delta(f))=2(d(f)+c(f))+\mu(\Sigma(f))$ where $d(f)$ is the
number of nodes and $c(f)$ is the number of cusps of $f$. As,
$\mu(\Sigma(f))=108$ and
$$\displaystyle{c(f)+d(f)=\dim_\mathbb{C}\frac{\mathcal{O}_{\mathbb{C}^2,0}}{ {F}_1(f)}}=2886,$$
we have that $\mu(\Delta(f))=5880$. In this case, by computation, we
obtain $c(f)=147$ and, therefore, the number of ordinary double
points is $d(f)=2739$.

\subsection{Topological classification  in $\mathcal O_2^2$}\label{secTopoClass}

For corank 2 map germs from $\mathbb C^2$ to $\mathbb C^2$ Gaffney
and Mond in \cite{GM2} ask the following question:
\begin{center}
{\it How many different topological types are contained in a given
$\mathcal K(xy,x^a+y^b)$-orbit?}
\end{center}

Miranda, Saia and Soares in \cite{MSS} show that for all pairs $(a,b)$, excluding
$(2,3)$ and $(2,5)$, there  exist a non finite number of distinct topological types in each
$\mathcal K$-orbit, that is, for $f(x,y)=(xy, \alpha x^a+ \beta y^b)$ there is at
least one family such that each element in the family is $\mathcal A$-finitely
determined germ and any two of them are not $C_0-\mathcal A$-equivalent. For
$(a,b)=(2,3)$ there is only one, and for $(a,b)=(2,5)$ there are two distinct
topological orbits.

In order to illustrate the use of the previous algorithm we fix
$(a,b)=(3,4)$. Let
$f_{u,v,w}(x,y)=(xy,y^4+x(x+y)^2+uxy^3+vx^3y+wx^4)$ be a family of
such maps. We use our implementation to obtain a presentation matrix
of $f_{u,v,w}$ restricted to the singular set. Now, the discriminant
of this map is given by the $0$-th Fitting ideal, and we have:
\begin{itemize}
\item if $u \neq v-2w+2$, then the Milnor number of discriminant curve is $\mu(\Delta(f_{u,v,w}))=54$.
\item if $u=v-2w+2$, then the family is not $\mathcal{A}-$finitely determined.
\end{itemize}

We set $u=2, v=w=0$ and consider the one parameter family
$$f_{s}(x,y)=(xy,y^4+x(x+y)^2+2xy^3+x^sy^{s+1})$$ with $s>3$. Now $f_s$ is
${\mathcal K}$-equivalent to $(xy,x^3+y^4)$ for all $s$, and the set
of generators of $\mathcal O_2$ as $\mathcal
O_{(\Sigma(f_s),0)}$-module via $\tilde f_s^*$ is $\{1, x, x^2, y,
y^2, y^3, y^4\}$. Using our implementation we obtain the following
presentation matrix for $f_{s_*}\mathcal{O}_{\Sigma(f_s)}$:

\bigskip

{\footnotesize \begin{tiny}
$\left[ \begin {array}{ccccccc} Y&-\frac{4}{3}\,X&0&-\frac{4}{3}\,X-\frac{4}{3}\,{X}^{s}&-\frac{10}{3}\,X&0&-\frac{7}{3}\\
\noalign{\medskip}-\frac{4}{3}\,{X}^{s+1}-\frac{4}{3}\,{X}^{2}&Y&-\frac{4}{3}\,X&-\frac{10}{3}\,{X}^{2}&0&-\frac{7}{3}\,X&0\\
\noalign{\medskip}-\frac{10}{3}\,{X}^{3}-{\frac {37}{30}}\,XY&-\frac{4}{3}\,{X}^{s+1}+\frac{6}{5}\,{X}^{2}&Y&\frac{6}{5}\,{X}^{s+1}+\frac{6}{5}\,{X}^{2}&0&0&{\frac {11}{10}}\,X\\
\noalign{\medskip}-\frac{5}{2}\,{X}^{2}&0&-\frac{7}{4}\,X&Y&-\frac{3}{4}\,{X}^{s}-\frac{3}{4}\,X&-X&0\\
\noalign{\medskip}0&-\frac{7}{4}\,{X}^{2}&0&-\frac{5}{2}\,{X}^{2}&Y&-\frac{3}{4}\,{X}^{s}-\frac{3}{4}\,X&-X\\
\noalign{\medskip}-{\frac {37}{12}}\,{X}^{3}&0&-\frac{4}{3}\,{X}^{2}&\frac{1}{3}\,XY&-\frac{5}{2}\,{X}^{2}&Y+\frac{2}{3}\,{X}^{2}&-\frac{3}{4}\,{X}^{s}-\frac{3}{4}\,X\\
\noalign{\medskip}{\frac {15}{8}}\,{X}^{3}+{\frac{27}{160}}\,XY& \beta &X^2(\frac {21}{16}-{X}^{s-1})
&\alpha &X \left( {\frac {9}{16}}\,{X}^{s}+Y \right) &0&Y-{\frac {63}{160}}\,X\end {array} \right]
$
\end{tiny}}

\medskip

\noindent where $\alpha =
-\frac{1}{2}\,XY-\frac{9}{4}\,{X}^{3}-{\frac
{9}{40}}\,{X}^{2}-\frac{1}{40}\,{X}^{s} \left( 9\,X-10\,Y \right)$
and $\beta =-\frac{5}{2}\,{X}^{3}-{\frac
{9}{40}}\,{X}^{2}-\frac{3}{4}\,XY.$

\medskip

Computing the Fitting ideals we obtain the following topological
invariants: $$\mu(\Delta(f_s))=2s+50,\ \ \ d(f_s)=s+14,\ \ \
\mu(\Sigma(f_s))=4\ \ \ \mbox{and}\ \ \ c(f_s)=9.$$ Therefore, for each
$s$ we have a distinct topological type.


\subsection{Number of triple points}

Consider the corank 2 quasi-homogeneous map germ $f: (\C^3,0)
\longrightarrow (\C^3,0)$, given by $$(x,y,z)\longmapsto
(x,yz,z^{18}+y^2+x^4z),$$ with weights $(17,36,4)$ and degrees
$(17;40;72)$. It is easy to verify that the Jacobian ideal is
$J(f)=\langle 18z^{18}+x^4z-2y^2\rangle$ and a set of generators for
$\mathcal{O}_{3}/({\widetilde{f}+J(f)})$ is $\{ 1,y,z,...,z^{18}
\}$. The time to get the presentation matrix was about $44$ seconds.

If $I(A_{1,1})$ is the ideal that defines the ordinary double points
in the target, $\sharp A_{(1,2)}, \sharp A_{(1,1,1)}$ and $\sharp
A_{3}$ are, respectively, the number of transversal intersection
between cuspidal edges and ordinary planes, the number of ordinary
triple points and the number of swallowtail as in the classical
Arnold notation, then by \cite{JMS}[Proposition 4.6 and 4.9] we have
that
$$\sharp A_{(1,2)}+ \sharp A_{(1,1,1)}+ \sharp A_3=\dim_{\C} \frac{ {\mathcal O}_{ 3 }}  { {  F}_2(f)}= 7368\ \ \ \mbox{and}\ \ \ \
\sharp A_{(1,1,1)}=\dim_{\C} \frac{ {\mathcal O}_3} { (
I(A_{1,1})^2:{  F}_0(f))}=5120.$$

For the quasi-homogeneous case, Ohmoto in \cite{Toru} presents
formulae to compute these $0$-stable singularities from $(\C^3,0)$
to $(\C^3,0)$. In this case, we obtain
$$\sharp A_3=136; \ \ \ \sharp A_{1,1,1}= 5120; \ \ \  \sharp A_{1,2} = 2112,$$
which can be confirmed directly using the previous presentation
matrix.

\subsection{Target multiple points}\label{secTargetMultiplePoints}

A simple but important application of our algorithm is to compute
the multiple spaces in the target $M_k(f)$ of a finite map germ
$f\colon (\C^n,0)\to(\C^{n+1},0)$. As mentioned in the introduction,
$M_k(f)$ is the zero set of the ideal $F_{k-1}(f)$ in $\cO_{n+1}$.
\begin{example}\label{ExDF}
{\rm Let $f\colon (\C^2,0)\to (\C^3,0)$ be the corank two map germ
given by $$(x,y)\mapsto(x^2,y^2,x^3+y^3+xy).$$ Our implementation
yields the following presentation matrix of $f_*\cO_2$:
$$\Lambda_f = \left[ \begin {array}{cccc} Y&-X_{{1}}&-X_{{2}}&-1
\\ \noalign{\medskip}-X_{{2}}Y-{X_{{1}}}^{2}&Y+X_{{1}}X_{{
2}}&{X_{{2}}}^{2}-X_{{1}}&0\\ \noalign{\medskip}-{X_{{2}}}^{2}-X_{{1}}
Y&-X_{{2}}+{X_{{1}}}^{2}&Y+X_{{1}}X_{{2}}&0
\\ \noalign{\medskip}-X_{{1}}X_{{2}}&-{X_{{2}}}^{2}&-{X_{{1}}}^{2}&Y\end {array} \right]
$$
and we obtain the following Fitting ideals (see Figure \ref{figure2}):

\begin{itemize}
\item ${ F}_0(f)=\langle X_1^2X_2^2-2X_1X_2Y^2+Y^4-2X_1^4X_2-2X_1X_2^4-8X_1^2X_2^2Y-2X_1^3Y^2-2X_2^3Y^2+X_1^6-2X_1^3X_2^3+X_2^6\rangle,$
 whose zeros define the image of $f$.

\item  ${ F}_1(f)= \langle X_2^2+X_1Y, Y+X_1X_2, -X_2+X_1^2\rangle \cap  \langle X_1+X_2-Y, X_2^2-X_2Y+Y^2\rangle 
 \cap \langle X_2+Y, X_1+Y \rangle \cap  \langle -X_1+X_2^2, Y+X_1X_2, X_1^2+X_2Y\rangle$.
This ideal defines the double point space of $f$ in the target.

\item  ${ F}_2(f)= \langle X_1,X_2,Y \rangle$, which defines the triple point in the image of $f$. As the codimension of $F_2(f)$ is $1$, this indicates that $f$ has exactly one triple point collapsed in the origin.
\end{itemize}

}
\end{example}
In the same way that the $\C$-codimension of $F_2(f)$ in $\cO_3$ measures the number of triple points collapsed in the origin of a map germ $f:(\C^2,0)\to(\C^3,0)$, we can use the algorithm to compute the number of $(n+1)$-tuple points collapsed in the origin of a map germ $(\C^n,0)\to (\C^{n+1},0)$. This invariant is just the codimension of $F_k(f)$ as a $\C$-vector subspace of $\cO_{n+1}$.

\subsection{Source double points of map germs $f:(\mathbb C^n,0)\to(\mathbb C^{n+1},0)$}\label{secSourceDoublePoints}

Let $f:(\mathbb C^n,0)\to(\mathbb C^{n+1},0)$ be a germ of finite
and generically-one-to-one map.  Following Mond
\cite{MondSomeRemarks}, the (lifted) double point space of $f$ is
the space $D^2(f)$ given by the ideal \[I^2(f)=(f\times
f)^*I_{n+1}+R(\alpha),\]
 where $I_{n+1}$ is the ideal defining the
diagonal of $\mathbb C^{n+1}\times \mathbb C^{n+1}$, and $R(\alpha)$
is given by the minors of any matrix $\alpha$, whose entries
$\lambda_{ij}\in \mathcal O_{2n}$ satisfy
\[f_j(x)-f_j(x')=\sum_{i=1}^n\alpha_{ij}(x,x')(x-x'),\]
 for all $1\leq
j\leq n+1$. In \cite{JN} it is shown that $D^2(f)$ is a Cohen-Macaulay
space of dimension $n-1$ (an extension of this result for map germs $(\C^n,0)\to (\C^p,0)$, with $n\leq p$, can be found in \cite{NP}). Set theoretically, $D^2(f)$ is given by
the pairs $(x,x')\in \C^n\times\C^n$, such that $f(x)=f(x')$ and, if
$x=x'$, then $f$ is singular at $x$. In \cite{MararMond} the source
double point space $D(f)\subset \C^n$ is defined as the image of
$D^2(f)$ by the projection on the first component
$\pi:(\C^n\times\C^n,0)\to(\C^n,0)$, that is:
$$D(f)=V(F_0(\pi\vert_{D^2(f)})).$$

The set $D(f)$ plays an important role. For instance, for map germs
from $\mathbb C^2$ to $\mathbb C^3$, it  characterizes finite
determinancy. More precisely: a map germ $f:(\C^2,0)\to(\C^3,0)$ is
finitely $\mathcal A$-determined if and only if the Milnor number
$\mu(D(f))$ is finite \cite{MararMond, MararNunoPenafortCor2Curves}.

Computing $D(f)$ for a map germ $(\C^n,0)\to(\C^p,0)$ can be quite
involved but, in the  case $p=n+1$, since $\pi\colon{D^2(f)}\to\C^n$
is a map from a Cohen-Macaulay space of dimension $n-1$ to $\C^n$,
we can use our algorithm to do so.
\begin{example}{\rm Consider the corank $2$ map germ $f(x, y) = (x^2, y^2, x^3 + y^3 + xy)$ of Example \ref{ExDF}. The set of double points $D^2(f)$ is given by the ideal
$I^2(f)$, generated by
\begin{enumerate}
\item[]$(x+u)(y+v), (x+u)(2y^2+2yv+2v^2+x+u),$
\item[]$(2x^2+2xu+2u^2+y+v)(y+v),$
\item[]$ x^2-u^2,$
\item[]$y^2-v^2,$
\item[]$x^3+y^3+xy-u^3-v^3-uv$.
\end{enumerate}

Take the projection $\pi:(\C^2 \times \C^2,0) \longrightarrow
(\C^2,0)$ given by $\pi(x,y,u,v)=(x,y)$. A basis of the vector space
$\frac{\mathcal O_{D^2(f)}}{\pi^*{\mathfrak m}}$
 is given by $\{1,y,u,v,v^2,v^3\}$. Using our implementation, we find the
following presentation matrix for $\pi_*\mathcal O_{D^2(f)}$:
$$\Lambda_{\pi} = \left[ \begin {array}{cccccc}
Y                         &-1                         &0                                                        &0                            &0                                                &0\\
0                               &Y                    &0                                                        &0                            &-1                                               &0\\
2XY+{X}^{3}   &0                          &Y-{X}^{2}                                    &X-X^{2}Y &-X^{2}                                   &-2\\
0                               &X^{2}              &0                                                        &Y+X^{2}        &             1                                   &0\\
0                               &X^2Y       &0                                                        &X^{2}Y         &Y                                          &1\\
\frac{1}{2}\,X^{2}      &-\frac{1}{2}\,X^{3}&\frac{1}{2}\,X+\frac{1}{2}\,X^{2}Y   &-\frac{1}{2}\,XY &-\frac{1}{2}\,X +\frac{1}{2}\,X^{2}Y&Y+\frac{1}{2}\,X^{2}\\
\end {array}\right]
$$

We rename $X=x$ and $Y=y$, since the target of $\pi$ is the source
of $f$, and then we obtain the following (Figure \ref{figure2}):
\begin{enumerate}
\item The ideal ${  F}_0(\pi|_{D^2(f)})=\langle(x^3+y^3)(x+y^2)(y+x^2)\rangle$, which defines the source double point space $D(f)$.

\item The ideal ${  F}_1(\pi|_{D^2(f)})= \langle x^2,xy,y^2\rangle$. We may regard the double
points of $\pi_{D^2(f)}$ as triple points of $f$. The codimension of ${  F}_1(\pi|_{D^2(f)})$ is $3$, corresponding to the number of source points in an ordinary triple point, which here have collapsed at $0$.

\end{enumerate}

 }
\begin{figure}
\begin{center}
\includegraphics[scale=0.675]{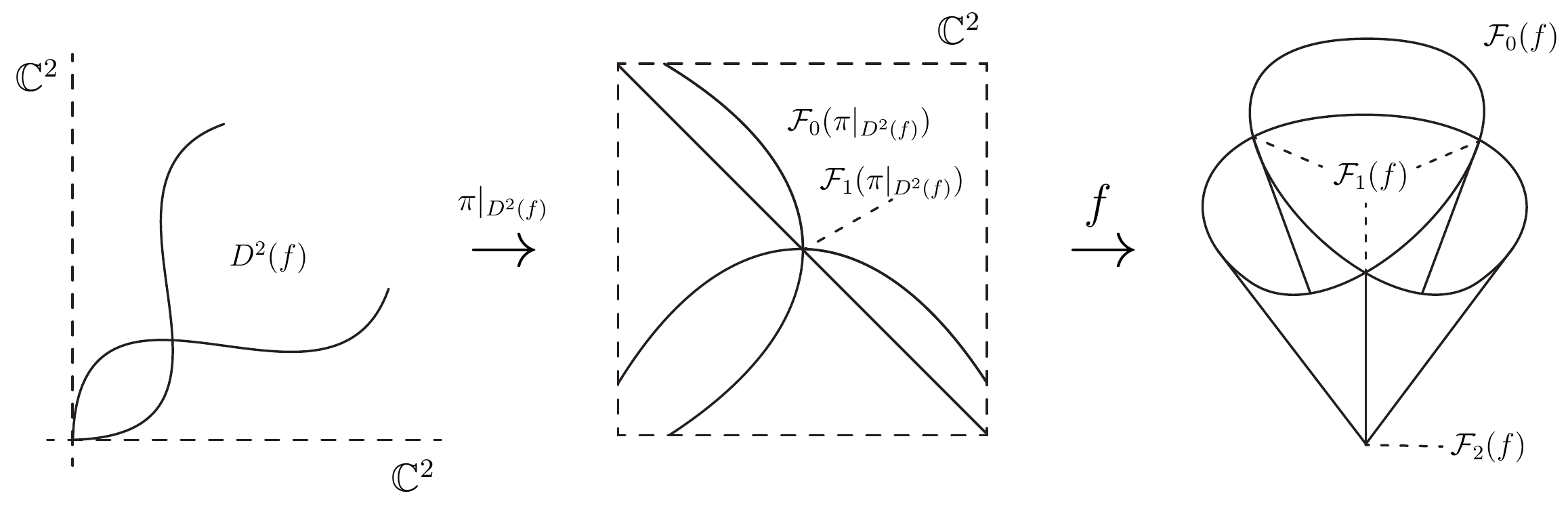}
\caption{\label{figure2}Lifted double points, double points and multiple points in the target.}
\end{center}
\end{figure}
 \end{example}

\begin{remark}
{\rm We can give two different analytic structures defining the $k$th source multiple point space of a map germ $f\colon (\C^n,0)\to (\C^p,0)$. In the first one, we regard the $k$th source multiple point space as $\pi (D^k(f))$, where $D^k(f)\subseteq  (\mathbb C^{n})^k$ is the lifted $k$th multiple point space, and $\pi$ is the projection on the first copy of $\mathbb C^n$. Thus, the defining ideal is $F_0(\pi_{D^k(f)})$. The second structure is given by defining the source multiple points as the preimage by $f$ of $M_k(f)$. It is an open problem to decide whether or not these analytic structures coincide. In the previous example, we computed ${  F}_0(\pi|_{D^2(f)})=\langle(x^3+y^3)(x+y^2)(y+x^2)\rangle$, which is precisely the preimage by $f$ of the ideal ${ F}_1(f)$ computed in Example \ref{ExDF}.}
\end{remark}
 {\bf Acknowledgements:}
 The authors are grateful to Juan Jos\'e Nu\~no Ballesteros, for many useful suggestions to the present work.

\vspace{1cm}

\begin{tabular}{lllll}
Hernandes, M. E. &  & Miranda, A. J. &  & Pe\~nafort Sanchis, G. \\
mehernandes$@$uem.br & & aldicio@ufu.br & & guillermo.penafort@uv.es \\
DMA-UEM & & FAMAT-UFU & & IMPA\\
Av. Colombo 5790 & & Av. João Naves de Ávila, 2121 & & Estr. Dona Castorina 110, 22460-320\\
Maringá-PR 87020-900 & & Uberlândia - MG 38408-100 & & Rio de Janeiro, \\
Brazil & & Brazil & & Brazil
\end{tabular}

\end{document}